\documentclass[11pt]{article}

\usepackage[a4paper]{anysize}\marginsize{2.5cm}{2.5cm}{1.3cm}{2cm}
\pdfpagewidth=\paperwidth \pdfpageheight=\paperheight
\usepackage{amsfonts,amssymb,amsthm,amsmath,eucal}
\usepackage{pgf}
\usepackage{bbm}
\usepackage{tikz} 
\usepackage{mathrsfs}
\usetikzlibrary{arrows}
\usepackage{color}
\usepackage{float}
\usepackage{subfigure}
\usepackage{caption}
\usepackage{graphicx}

\definecolor{qqqqff}{rgb}{0.,0.,0.}

\theoremstyle{plain}
\newtheorem{thm}{Theorem}[section]
\newtheorem{theorem}[thm]{Theorem}

\newtheorem{lemma}[thm]{Lemma}
\newtheorem{proposition}[thm]{Proposition}

\newtheorem{conjecture}[thm]{Conjecture}

\theoremstyle{definition}
\newtheorem{definition}[thm]{Definition}
\newtheorem{remark}[thm]{Remark}
\newtheorem{example}[thm]{Example}

\newtheorem{thevarthm}[thm]{\varthmname}

\newenvironment{varthm*}[1]{\trivlist\item[]{\bf #1.}\it}{\endtrivlist}

\renewcommand\geq{\geqslant}

\renewcommand\leq{\leqslant}

\newcommand\be{\begin{eqnarray*}}
\newcommand\ee{\end{eqnarray*}}

\newcommand\Q{\mathbb Q}
\newcommand\R{\mathbb R}
\newcommand\C{\mathbb C}

\renewcommand\P{\mathbb P}

\newcommand\calo{{\mathcal O}}
\newcommand\cali{{\mathcal I}}

\newcommand\newop[2]{\def#1{\mathop{\rm #2}\nolimits}}
\newop\log{log}
\newop\ord{ord}
\newop\Gal{Gal}
\newop\SL{SL}
\newop\Bl{Bl}
\newop\mult{mult}
\newop\mass{mass}
\newop\div{div}
\newop\codim{codim}
\newop\sing{sing}
\newop\Zeroes{Zeroes}

\newcommand\alphahat{\widehat{\alpha}}
\newcommand\eps{\varepsilon}

\def\keywordname{{\bfseries Keywords}}%
\def\keywords#1{\par\addvspace\medskipamount{\rightskip=0pt plus1cm
\def\and{\ifhmode\unskip\nobreak\fi\ $\cdot$
}\noindent\keywordname\enspace\ignorespaces#1\par}}
\def\subclassname{{\bfseries Mathematics Subject Classification
(2000)}\enspace}
\def\subclass#1{\par\addvspace\medskipamount{\rightskip=0pt plus1cm
\def\and{\ifhmode\unskip\nobreak\fi\ $\cdot$
}\noindent\subclassname\ignorespaces#1\par}}

\def\endproof{\hspace*{\fill}\endproofsymbol\endtrivlist}

\def\endproofsymbol{\frame{\rule[0pt]{0pt}{6pt}\rule[0pt]{6pt}{0pt}}}

\begin{document}

\author{M.~Dumnicki, T.~Szemberg, J.~Szpond}
\title{Local effectivity in projective spaces}
\date{\today}
\maketitle
\thispagestyle{empty}

\begin{abstract}
   In this note we introduce a Waldschmidt decomposition of divisors which might be
   viewed as a generalization of Zariski decomposition based on the effectivity rather than
   the nefness of divisors. As an immediate application we prove a recursive formula providing
   new effective lower bounds on Waldschmidt constants
   of very general points in projective spaces. We use these bounds in order to verify
   Demailly's conjecture in a number of new cases.
\keywords{Chudnovsky conjecture, Demailly conjecture, Waldschmidt constants}
\subclass{MSC 14C20 \and MSC 14J26 \and MSC 14N20 \and MSC 13A15 \and MSC 13F20}
\end{abstract}


\section{Introduction}
   Let $X$ be a smooth projective variety and let $L$ be an ample line bundle
   on $X$. The concept of the \emph{local positivity} of $L$ has been coined
   by Demailly, who introduced in \cite{Dem90} the following invariants
   measuring in effect the local positivity.
\begin{definition}[Seshadri constant]
   Let $X$ be a smooth projective variety and let $L$ be an ample line bundle
   on $X$. Let $P\in X$ be a fixed point and let $f:\Bl_PX\to X$ be the blow
   up of $X$ at $P$ with the exceptional divisor $E$. The real number
   $$\eps(X;L,P)=\sup\left\{t\in\R:\;f^*L-tE\;\mbox{ is nef}\right\}$$
   is the \emph{Seshadri constant} of $L$ at $P$.
\end{definition}
   Thus $\eps(X;L,P)$ is the value of $t$ for which the ray $f^*L-tE$
   hits the boundary of the nef cone on $\Bl_P X$.
   It is natural to introduce a similar invariant, which gives the value
   of $t$, where the ray $f^*L-tE$ hits the boundary of the pseudo-effective cone on $\Bl_P X$.
   We consider this invariant (more precisely its reciprocal introduced in Definition \ref{def:Waldschmidt constant})
   as a way to measure the \emph{local effectivity}
   of $L$.
\begin{definition}[The $\mu$-invariant]
   Let $X$ be a smooth projective variety and let $L$ be an ample line bundle
   on $X$. Let $P\in X$ be a fixed point and let $f:\Bl_PX\to X$ be the blow
   up of $X$ at $P$ with the exceptional divisor $E$. The real number
   $$\mu(X;L,P)=\sup\left\{t\in\R:\;f^*L-tE\;\mbox{ is effective}\right\}$$
   is the \emph{$\mu$-invariant} of $L$ at $P$.
\end{definition}
   Both notions can be easily generalized replacing the point $P$ by an
   arbitrary subscheme $Z\subset X$ and taking $f:\Bl_ZX\to X$ to be the blow up
   of $X$ along the ideal sheaf $\cali_Z\subset\calo_X$. We denote the exceptional
   divisor of $f$ again by $E$.

   Whereas the $\mu$-invariant $\mu(X;L,Z)$ is not much present in the literature,
   its reciprocal is the well-known Waldschmidt constant of $Z$.
   We define first the \emph{initial degree of $Z$ with respect to $L$} as
   $$\alpha(X;L,Z)=\min\left\{d:\; df^*L-E\;\mbox{ is effective}\right\}.$$
   For an integer $m\geq 1$, let $mZ$ denote the subscheme defined by the symbolic power
   $\cali_Z^{(m)}$ of $\cali_Z$, see \cite[Definition 9.3.4]{PAG}.
   Then the asymptotic version of the initial degree is the following.
\begin{definition}[Waldschmidt constant]\label{def:Waldschmidt constant}
   Let $X$ be a smooth projective variety and let $L$ be an ample line bundle on $X$.
   Let $Z\subset X$ be a subscheme. The real number
   $$\alphahat(X;L,Z)=\inf_{m\geq 1}\frac{\alpha(X;L,mZ)}{m}$$
   is the \emph{Waldschmidt constant} of $Z$ with respect to $L$.
\end{definition}
\begin{remark}\rm
   Since the numbers $\alpha(X;L,mZ)$ for $m\geq 1$ form a subadditive sequence, i.e.
   there is
   $$\alpha(X;L,(k+\ell)Z)\leq \alpha(X;L,kZ) + \alpha(X;L,\ell Z)$$
   for all $k$ and $\ell$, the infimum in Definition\ref{def:Waldschmidt constant}
   exists and moreover we have
   $$\alphahat(X;L,Z)=\lim_{m\to\infty}\frac{\alpha(X;L,mZ)}{m}.$$
\end{remark}
   Waldschmidt constants appear in different guises in various branches
   of mathematics. Apparently, they were first considered in complex analysis
   in connection with estimates on the growth order of holomorphic functions,
   see \cite{Wal77}. In this setup $X$ is simply $\C^n$ or $\P^n$. We prefer
   the homogeneous approach here. Then the polarization $L$ is just the
   hyperplane bundle $\calo_{\P^N}(1)$. Let $I$ be a non-zero, proper homogeneous ideal
   in the polynomial ring $\C[x_0,\ldots,x_N]$. The \emph{initial degree}
   of $I$ is
   $$\alpha(\P^N;I)=\min\left\{d:\; (I)_d\neq 0\right\},$$
   where $(I)_d$ denotes the degree $d$ part of $I$.
   The Waldschmidt constant of $I \subset \C[x_0,\ldots,x_N]$ is then
   $$\alphahat(\P^N;I)=\inf_{m\geq 1}\frac{\alpha(\P^N;I^{(m)})}{m},$$
   which of course agrees with Definition \ref{def:Waldschmidt constant}.
   In recent years there has been considerable interest in Waldschmidt constants in general,
   see e.g. \cite{DHST14}, \cite{BCGHJNSvTT16}, \cite{MosHag16}, \cite{FGHLMS17}.
   Special attention has been given to the following
   Conjecture stated originally by Demailly in \cite[p. 101]{Dem82}.
   It has been formulated recently by Harbourne and Huneke in \cite[Question 4.2.1]{HaHu13}.
   Apparently the authors were not aware of Demailly's work. We use again the projective version.
\begin{conjecture}[Demailly]\label{conj:Demailly}
   Let $Z\subset\P^N$ be a finite set of points and let $I$ be the homogeneous
   saturated ideal defining $Z$. Then for all $m\geq 1$
   \begin{equation}\label{eq:Demailly Conjecture}
      \alphahat(\P^N;I)\geq \frac{\alpha(\P^N;I^{(m)})+N-1}{m+N-1}.
   \end{equation}
\end{conjecture}
   For $m=1$ the Conjecture of Demailly reduces to the statement which is best known as the Conjecture of Chudnovsky,
   see \cite[Problem 1]{Chu81}, to the effect that the inequality
   \begin{equation}\label{eq:Chudnovsky Conjecture}
      \alphahat(\P^N;I)\geq \frac{\alpha(\P^N;I)+N-1}{N}.
   \end{equation}
   holds for all ideals defining finite sets of points in $\P^N$.
   Demailly's Conjecture for $\P^2$ has been proved by Esnault and Viehweg
   using methods of complex projective geometry, see \cite[In\'egalit\'e A]{EsnVie83}.

   In the present note, we provide lower bounds on Waldschmidt constants
   of sets of general points in projective spaces and obtain as a corollary
   a proof of the Demailly's Conjecture in certain cases, see Theorem \ref{thm:DC ok for geq mN}.
   The new tool developed in this note is the concept of Waldschmidt decomposition
   introduced in Section \ref{sec:WD}. Our main results are
   Theorem \ref{thm:main} which gives an iterative way to control
   Waldschmidt constants of very general points and Proposition \ref{prop:distribution k s-k}
   which is an effective criterion derived from Theorem \ref{thm:main}.

\paragraph{Convention and notation.}
   We work throughout over the field $\C$ of complex numbers.

\section{Waldschmidt decomposition}\label{sec:WD}
   The numerical meaning of the Waldschmidt constant $\alphahat(X;L,Z)$
   is that if $D\in|kL|$ is an effective divisor vanishing along $Z$
   with multiplicity $m$, then
   $$\frac{k}{m}\geq \alphahat(X;L,Z).$$
   This condition extends easily to effective $\R$-divisors. Indeed, let
   $D=\sum\delta_iD_i$ be an effective $\R$-divisor with $D\equiv\delta L$
   for some $\delta>0$. Then $\mult_ZD=\sum\delta_i\mult_ZD_i$ and
   $$\frac{\delta}{\mult_ZD}\geq \alphahat(X;L,Z).$$
   In this section we introduce certain decomposition of a divisor,
   depending on its numerical properties. We call it the Waldschmidt
   decomposition as it is governed by Waldschmidt constants.
   This decomposition can be viewed as a higher dimensional version
   of the Bezout decomposition defined in \cite[Section 2.1]{DST16}.
   Whereas it is possible to define it on arbitrary varieties, we restrict
   our approach here to $\P^N$ and its linear subspaces. In this setting
   the definition is most transparent.
\begin{definition}[Waldschmidt decomposition in $\P^N$]\label{def:Waldschmidt decomposition}
   Let $H\cong \P^{N-1}$ be a hyperplane in $\P^N$ and let $Z$
   be a subscheme in $H$. Let $D$ be a divisor of degree $d$ in $\P^N$.
   The \emph{Waldschmidt decomposition of $D$ with respect to $H$ and $Z$} is the sum
   of $\R$-divisors
   $$D=D'+\lambda\cdot H$$
   such that $\deg(D')=d-\lambda,$
   \begin{equation}\label{eq:Waldschmidt decomposition cond}
      \frac{d-\lambda}{\mult_ZD'}\geq\alphahat(H;\calo_H(1),Z)
   \end{equation}
   and $\lambda$ is the least non-negative real number such that \eqref{eq:Waldschmidt decomposition cond}
   is satisfied.
\end{definition}
   Of course, it may happen that $\lambda=0$ in Definition \ref{def:Waldschmidt decomposition}.
   This number is positive, if the restriction of $D$ to $H$ would produce a divisor in $|\calo_H(1)|$
   violating the inequality \eqref{eq:Waldschmidt decomposition cond}. Thus $\lambda$ is the least multiplicity such that
   $H$ is numerically forced to be contained in $D$ with this multiplicity. It may well happen
   that the divisor $D'$ still contains $H$ as a component.
\begin{remark}
   The definition of the Waldschmidt decomposition with respect to $H$ can be extended to a finite number of hyperplanes $H_1,\ldots,H_s$.
\end{remark}

\section{The main result}
   In this section we state our main result. The statement is motivated by the proof of the following lower
   bound on Waldschmidt constants presented in \cite[Theorem 3]{DT16}.
\begin{theorem}[Lower bound on Waldschmidt constants]\label{thm:lower bound}
   Let $I$ be the saturated ideal of a set of $r$ very general points in $\P^N$. Then
   $$ \alphahat(\P^N;I)\geq\lfloor \sqrt[N]{r}\rfloor.$$
\end{theorem}
   It is expected that for $r$ sufficiently big, there is actually the equality $\alphahat(\P^N;I)=\sqrt[N]{r}$
   but this statement seems out of reach with present methods.
\begin{theorem}\label{thm:main}
   Let $H_1,\ldots,H_s$ be $s\geq 2$ mutually distinct hyperplanes in $\P^N$.
   Let $a_1,\ldots,a_s \geq 1$ be real numbers such that
   \begin{equation}\label{eq:inequalities on a_i}
   1 - \sum_{j=1}^{s-1} \frac{1}{a_j} > 0
   \end{equation}
   and
   \begin{equation}\label{eq:inequality on a_s}
   1 - \sum_{j=1}^{s} \frac{1}{a_j} \leq 0.\\
   \end{equation}
   Let
   $$Z_i=\left\{P_{i,1},\ldots,P_{i,r_i}\right\}\in H_i\setminus \bigcup_{j\neq i}H_j$$
   be the set of $r_i$ points
   such that
   \begin{equation}\label{eq:bound on alpha H_i}
      \alphahat(H_i;Z_i)\geq a_i
   \end{equation}
   and let $Z=\bigcup_{i=1}^s Z_i$.
   Finally, let
   \begin{equation}\label{eq:q}
      q:=\left(1-\sum\limits_{j=1}^{s-1} \frac{1}{a_j}\right) \cdot a_s+s-1.
   \end{equation}
   Then
   $$\alphahat(\P^N;Z)\geq q.$$
\end{theorem}
\proof
First observe that, for any $t=1,\dots,s-1$, by \eqref{eq:inequalities on a_i} we have
$$1 - \sum_{j=1}^{t} \frac{1}{a_j} > 0.$$
Multiplying by $a_t$, moving $a_t/a_t=1$ to the right hand side and making some preparation we get
$$a_t - \sum_{j=1}^{t-1} \frac{a_t}{a_j} > \left(1 - \sum_{j=1}^{t-1} \frac{1}{a_j} \right) + \sum_{j=1}^{t-1} \frac{1}{a_j}.$$
Dividing by $1-\sum\limits_{j=1}^{t-1} \frac{1}{a_j}$ we get
\begin{equation}\label{eq:incmore}
a_t > 1 + \frac{\sum\limits_{j=1}^{t-1} \frac{1}{a_j}}{1 - \sum\limits_{j=1}^{t-1}\frac{1}{a_j}}
\end{equation}
for $t \leq s-1$. Similarly, starting with \eqref{eq:inequalities on a_i}, we get
\begin{equation}\label{eq:incmore2}
a_s \leq 1 + \frac{\sum\limits_{j=1}^{s-1} \frac{1}{a_j}}{1 - \sum\limits_{j=1}^{s-1}\frac{1}{a_j}}.
\end{equation}

   We assume to the contrary that there is a divisor $D$ of degree $d$ in $\P^N$ vanishing
   to order at least $m$ at all points of $Z$ such that
   \begin{equation}\label{eq:p lower than q}
      p:=\frac{d}{m} < q.
   \end{equation}
   It is convenient to work with the $\Q$-divisor $\Gamma=\frac1m D$, which is of degree $p$ and has multiplicities
   at least $1$ at every point of $Z$.

\medskip
\noindent
   \textbf{Step 0.}\\
   Let $\Gamma=\Gamma'+\sum_{i=1}^s\lambda_i H_i$ be the Waldschmidt decomposition of $\Gamma$ with respect
   to $H_1,\ldots,H_s$ and $Z_1,\ldots,Z_s$ respectively. The conditions \eqref{eq:Waldschmidt decomposition cond} and \eqref{eq:bound on alpha H_i}
   imply then that
   \begin{equation}\label{eq:Bezout condition for p}
   \left\{\begin{array}{ccrcl}
      (\ref{eq:Bezout condition for p}.1) && p-\sum\limits_{i=1}^s\lambda_i & \geq & a_1(1-\lambda_1)\\
      (\ref{eq:Bezout condition for p}.2) && p-\sum\limits_{i=1}^s\lambda_i & \geq & a_2(1-\lambda_2)\\
      \vdots && \vdots &&\\
      (\ref{eq:Bezout condition for p}.s) && p-\sum\limits_{i=1}^s\lambda_i & \geq & a_s(1-\lambda_s)\\
   \end{array}\right.
   \end{equation}
   We will show that the conditions in \eqref{eq:inequalities on a_i}, \eqref{eq:inequality on a_s}, \eqref{eq:p lower than q}
   and \eqref{eq:Bezout condition for p} cannot hold simultaneously. This will provide the desired
   contradiction to the existence of $D$. The idea is first to achieve equalities in \eqref{eq:Bezout condition for p}.

\medskip
\noindent
   \textbf{Step 1.}\\
   Our first claim is that there exists $\lambda_1'\leq \lambda_1$ such that
   \begin{equation}\label{eq:Bezout condition for p 1}
   \left\{\begin{array}{ccrcl}
      (\ref{eq:Bezout condition for p 1}.1) && p-\lambda_1'-\sum\limits_{i=2}^s\lambda_i & = & a_1(1-\lambda_1')\\
      (\ref{eq:Bezout condition for p 1}.2) && p-\lambda_1'-\sum\limits_{i=2}^s\lambda_i & \geq & a_2(1-\lambda_2)\\
      \vdots &&\\
      (\ref{eq:Bezout condition for p 1}.s) && p-\lambda_1'-\sum\limits_{i=2}^s\lambda_i & \geq & a_s(1-\lambda_s)\\
   \end{array}\right.
   \end{equation}
   Indeed, we have
   $$p-\lambda_1-\sum_{i=2}^s\lambda_i  \geq  a_1(1-\lambda_1)$$
   from (\ref{eq:Bezout condition for p}.1).
   Decreasing $\lambda_1$ by $\eps$, the left hand side increases
   by $\eps$ as well, whereas the right hand side increases
   by $a_1\cdot\eps$. Since $a_1>1$ by (\ref{eq:inequalities on a_i}.1),
   there must exist $\eps\geq 0$ such that
   $$p-(\lambda_1-\eps)-\sum\limits_{i=2}^s\lambda_i  \;=\;  a_1(1-(\lambda_1-\eps)).$$
   We put $\lambda_1'=\lambda_1-\eps$. Note also that decreasing $\lambda_1$
   preserves the inequalities with indices $j=2,\ldots,s$ in \eqref{eq:Bezout condition for p}
   because the left hand sides of all these inequalities increase, while
   the right hand sides remain unaltered.

   In order to alleviate the notation, we drop the prime index by the new $\lambda_1$.

\medskip
\noindent
   \textbf{Step t (the induction step).}
   In the second step we assume that we found new $\lambda_1,\dots,\lambda_{t-1}$ such that the following holds:
   \begin{equation}\label{eq:step2a}
   \left\{\begin{array}{rcl}
      p-\sum\limits_{i=1}^{t-1} \lambda_i-\lambda_t-\sum\limits_{i=t+1}^s\lambda_i & = & a_1(1-\lambda_1)\\
      \vdots &&\\
      p-\sum\limits_{i=1}^{t-1} \lambda_i-\lambda_t-\sum\limits_{i=t+1}^s\lambda_i & = & a_{t-1}(1-\lambda_{t-1})\\
      p-\sum\limits_{i=1}^{t-1} \lambda_i-\lambda_t-\sum\limits_{i=t+1}^s\lambda_i & \geq & a_{t}(1-\lambda_{t})\\
      \vdots &&\\
      p-\sum\limits_{i=1}^{t-1} \lambda_i-\lambda_t-\sum\limits_{i=t+1}^s\lambda_i & \geq & a_{s}(1-\lambda_{s})\\
   \end{array}\right.
   \end{equation}
   Our aim is to push this one step further, to the situation, where (for new $\lambda_1,\dots,\lambda_t$) we will have
   at least $t$ equalities.

   Let
   $$C := p - \lambda_t - \sum_{i=t+1}^{s} \lambda_i.$$
   Solve the following system of equalities with respect to  $\lambda_1,\dots,\lambda_{t-1}$
   and a parameter $\lambda_t$.
   \begin{equation}\label{eq:step2b}
   \left\{\begin{array}{rcl}
      C-\sum\limits_{i=1}^{t-1} \lambda_i & = & a_1(1-\lambda_1)\\
      \vdots &&\\
      C-\sum\limits_{i=1}^{t-1} \lambda_i & = & a_{t-1}(1-\lambda_{t-1})\\
   \end{array}\right.
   \end{equation}
   Let $\lambda_1',\ldots,\lambda_{t-1}'$ be unique (by Lemma \ref{solutionlem}) solutions to that system.
   Again, by Lemma \ref{solutionlem},
   \begin{equation}\label{eq:step2c}
   \sum\limits_{i=1}^{t-1} \lambda_{i}' = \frac{ C\left(\sum\limits_{j=1}^{t-1} \frac{1}{a_j} \right) - (t-1) }{ \left(\sum\limits_{j=1}^{t-1} \frac{1}{a_j} \right) - 1}.
   \end{equation}
   Since $\lambda_t$ is hidden in $C$ (as $-\lambda_t$), decreasing $\lambda_t$ by $\varepsilon$ increases
   $\sum\limits_{j=1}^{t-1} \lambda_i'$ by
$$ \varepsilon \left(  \frac{\sum\limits_{j=1}^{t-1} \frac{1}{a_j}}{\sum\limits_{j=1}^{t-1} \frac{1}{a_j}-1} \right).$$
   Thus the left hand side of the inequality \eqref{eq:step2a}.t increases by
   $$\varepsilon \left(1 + \frac{\sum\limits_{j=1}^{t-1} \frac{1}{a_j}}{1 - \sum\limits_{j=1}^{t-1}\frac{1}{a_j}}\right),$$
   which by \eqref{eq:incmore} is strictly less than $\varepsilon a_t$. In effect, decreasing $\lambda_t$, solving \eqref{eq:step2b} for $\lambda_1,\dots,\lambda_{t-1}$ gives a new sequence $\lambda_1',\dots,\lambda_t'$, with
\begin{itemize}
\item
preserved equalities $\eqref{eq:step2a}.1$ --- $\eqref{eq:step2a}.(t-1)$,
\item
left hand side of $\eqref{eq:step2a}.t$ increasing faster than the right hand side,
\item
left hand sides of $\eqref{eq:step2a}.(t+1)$ --- $\eqref{eq:step2a}.s$ increasing, while right hand sides remain unaltered.
\end{itemize}
As in Step 1, this suffices to obtain new $\lambda_1,\dots,\lambda_t$ with one more equality in \eqref{eq:step2a}.

\medskip
\noindent
   \textbf{Step s (the final step).}\\
   Assume that we have now $s-1$ equalities in \eqref{eq:step2a}, with the last inequality not necessaritly being an equality. We begin exactly as in the previous step. The only difference is that, by \eqref{eq:incmore2}, decreasing $\lambda_t$
forces the left hand side of the last inequality \eqref{eq:step2a} to increase \emph{faster} than the right hand side.
   Thus we may decrease $\lambda_t$ (altering $\lambda_1,\dots,\lambda_{t-1}$ to preserve equalities) to zero to obtain
   \begin{equation}\label{eq:Bezout condition for p s}
   \left\{\begin{array}{ccrcl}
      (\ref{eq:Bezout condition for p s}.1) && p-\sum\limits_{i=1}^{s-1}\lambda_i & = & a_1(1-\lambda_1)\\
      (\ref{eq:Bezout condition for p s}.2) && p-\sum\limits_{i=1}^{s-1}\lambda_i & = & a_2(1-\lambda_2)\\
      \vdots &&\vdots &&\\
      (\ref{eq:Bezout condition for p s}.(s-1)) && p-\sum\limits_{i=1}^{s-1}\lambda_i & = & a_{s-1}(1-\lambda_{s-1})\\
      (\ref{eq:Bezout condition for p s}.s) && p-\sum\limits_{i=1}^{s-1}\lambda_i & \geq & a_s\\
   \end{array}\right.
   \end{equation}
   It follows from Lemma \ref{solutionlem} that now
   \begin{equation}\label{eq:sum of lambda solutions}
      \sum\limits_{i=1}^{s-1}\lambda_i=\frac{pR-(s-1)}{R-1},
   \end{equation}
   where $R=\sum\limits_{j=1}^{s-1} 1/a_j$.
   From \eqref{eq:q} we have
   \begin{equation}\label{eq:loc3}
      q=(1-R)a_s+(s-1).
   \end{equation}
   Taking (\ref{eq:Bezout condition for p s}.s) into account we get
   $$q\leq (1-R)\left(p-\frac{pR-(s-1)}{R-1}\right)+(s-1)=p-Rp+pR-(s-1)+(s-1)=p.$$
   This contradicts however clearly \eqref{eq:p lower than q} and we are done.
\endproof

\section{Applications}

We will focus on Waldschmidt constants of sets of very general points in $\P^N$. The notation
$$\alphahat(\P^N;r)$$
denotes the Waldschmidt constant $\alphahat(\P^N;I)$ of a radical ideal $I$ of $r$ very general points in $\P^N$.

\begin{theorem}
\label{stepback}
Let $N \geq 2$, let $k \geq 1$ be an integer. Assume that for some integers
$r_1,\dots,r_{k+1}$ and rational numbers $a_1,\dots,a_{k+1}$ we have
$$\alphahat(\P^{N-1};r_j) \geq a_j \text{ for } j=1,\dots,k+1,$$
$$k \leq a_j \leq k+1 \text{ for } j=1,\dots,k, \quad a_1 > k, \quad a_{k+1} \leq k+1.$$
Then
$$\alphahat(\P^N;r_1+\ldots+r_{k+1}) \geq \left( 1-\sum\limits_{j=1}^{k} \frac{1}{a_j} \right) a_{k+1}+k.$$
\end{theorem}

\proof
We combine Theorem \ref{thm:main} and the specialization. We take hyperplanes $H_1,\dots,H_{k+1}$ and specialize
$r_j$ points to a set $Z_j \subset H_j$ for $j=1,\dots,k+1$, so that the points in $Z_j$ are in very general position on $H_j$. Hence
$$\alphahat(H_j;Z_j) = \alphahat(\P^{N-1};r_j).$$
To check that (\ref{eq:inequalities on a_i}) is satisfied, we compute
$$\sum_{j=1}^{k} \frac{1}{a_j} < \sum_{j=1}^{k} \frac{1}{k} = 1$$
since $a_j \geq k$ and $a_1 > k$.
Similarly we check that (\ref{eq:inequality on a_s}) holds,
$$\sum_{j=1}^{k+1} \frac{1}{a_j} \geq \sum_{j=1}^{k+1} \frac{1}{k+1} = 1.$$
The inequalities (\ref{eq:bound on alpha H_i}) are satisfied by assumptions. Thus the Waldschmidt constant of
specialized points is bounded as desired, hence for points in the very general position the bound also holds.
\endproof

\begin{example}
We bound from below $\alphahat(\P^3;20)$. Let $k=2$ (in fact, it is very easy to find the suitable $k$ in general; it must satisfy
$k^N < r < (k+1)^N$, where $r$ is the number of points in $\P^N$). Then we look for integers $r_1$, $r_2$ and $r_3$ and rational numbers
$a_1$, $a_2$, $a_3$ satisfying the assumptions of Theorem \ref{stepback}. Since we want to bound $\alphahat(\P^3;20)$, it must be
$$r_1 + r_2 + r_3 \leq 20.$$
Since $\alphahat(\P^2;r_1) \geq a_1 > 2$, we see that $r_1 > 4$. Similarly $r_2 \geq 4$, $r_3 \geq 1$.
Moreover, from $a_j \leq 3$ we see that we may restrict ourselves to the case when $r_j \leq 9$. Since $\alphahat(\P^2;r)$
is known for $r \leq 9$, it suffices to search through all the possibilities $(r_1,r_2,r_3)$, compute $(a_1,a_2,a_3)$ for each of them
and get the bound. This can be done by hand in principle. We have used a simple computer program to do the dully
calculations for us. As a result, for
$$r_1=8, \qquad r_2=8, \qquad r_3=4$$
we get
$$a_1 = 48/17, \qquad a_2 = 48/17, \qquad a_3=2.$$
   Thus, from the formula, $\alphahat(\P^3;20) \geq 31/12 \simeq 2.583$. Note that the upper bound is $\sqrt[3]{20} \simeq 2.714$.
\end{example}

\subsection{A recursive approach}
Now we study a much harder example which allows us to discuss some algorithmic issues.

\begin{example}
We want to bound $\alphahat(\P^4;180)$. Since now $N=4$, we get immediately $k=3$, since then $k^N < 180 < (k+1)^N$.
We are interested in sequences of integers
$$(r_1,r_2,r_3,r_4) \text{ with } r_1+r_2+r_3+r_4 \leq 180.$$
As before, we have additional constraints. Since $\alphahat(\P^{N-1};r_j) \geq a_j \geq k$, we get (in general) that
$r_j \geq k^{N-1}$. In our situation this gives $r_2,r_3 \geq 27$, $r_1 \geq 28$, $r_4 \geq 1$. It is reasonable to restrict
to $r_j \leq (k+1)^{N-1}$, so in our case, $r_j \leq 64$.

The first problem we encounter here is the number of sequences $(r_1,r_2,r_3,r_4)$ with above properties. But this can be
(in the case studied here, $N=4$, $r=180$) easily managed by a suitable computer program.
   What requires much more attention is
   coming up with good bounds $a_j$ for $\alphahat(\P^3;r_j)$. These constants are not known, except for several cases: $2$ for $\alphahat(\P^3;8)$,
$3$ for $\alphahat(\P^3;27)$ and $4$ for $\alphahat(\P^3;64)$. So the first approach is to use only numbers $r_j$
of the form $\ell^{N-1}$, which is weak, but manageable (we will address this later, in Proposition \ref{prop:distribution k s-k}).
Taking
$$r_1 = 64, \qquad r_2 = 64, \qquad r_3 = 27, \qquad r_4 = 8$$
we get
$$a_1 = 4, \quad a_2 = 4, \quad a_3 = 3, \quad a_4 = 2, \quad \text{thus} \quad \alphahat(\P^4;180) \geq \frac{10}{3} \simeq 3.333.$$

Using again a computer program we can find, as in the previous example, all necessary bounds for $\alphahat(\P^{N-1};\widetilde{r})$
for $\widetilde{r} = 1,\dots,(k+1)^{N-1}$. In our case it requires $64$ computations to find a bound in $\P^3$. Each of them requires
again looking for sequences satisfying certain properties and then going down to $\P^2$. In effect, the run time grows exponentially
when $N$ is increased. For $\alphahat(\P^2;\widetilde{r})$, however, a much better idea is to use known best bounds, e.g.,
   \cite[Theorem 2.2 and discussion thereafter]{HarRoe09}.

   Coming back to our case, with the help of a computer program, which run several minutes, all possibilities were scanned
   and the best results were found taking
$$r_1 = 52, \qquad r_2 = 52, \qquad r_3 = 49, \qquad r_4 = 27.$$
   Again with a computer we obtain
$$a_1 = a_2 = \frac{17457}{4816}, \quad a_3 = \frac{63495}{17974}, \quad a_4 = 3, \quad \text{thus} \quad \alphahat(\P^4;180) \geq 3.495.$$
In fact, the last number is exactly $430502824/123159135$. Observe that the upper bound is $\sqrt[4]{180} \simeq 3.663$.

From the above considerations we conclude that checking all partitions of $r$ into $k+1$ numbers would
take too much time for bigger $N$. To make this faster and manageable even in the case, e.g., $N \geq 100$ we must drastically reduce
the number of subcases. The radical idea is to consider only one distribution, and go down to $\P^{N-1}$ with only one case.

Observe that we look for the numbers $a_1,\dots,a_{k+1}$ such that
$$\left( 1 - \sum_{j=1}^{k} \frac{1}{a_j} \right) a_{k+1}$$
is as big as possible. The numbers $a_j$ are good bounds for $\alphahat(\P^{N-1};r_j)$, so we may as well assume, that
they are close to $\sqrt[N-1]{r_j}$ or even pretend they are equal.

   We consider first the expression
\begin{equation}\label{fracterm}
\sum_{j=1}^{k} \frac{1}{\sqrt[N-1]{r_j}}.
\end{equation}
   For all partitions $r_1+\dots+r_k = const$, we want \eqref{fracterm} to be as small as possible. Without
going into details, this forces all numbers $r_j$ to be nearly equal. Therefore we want to maximize
$$\left(1-\frac{k}{\sqrt[N-1]{r_1}}\right) \sqrt[N-1]{r_{k+1}}$$
   under the condition
$$kr_1 + r_{k+1} = r,$$
or, which is much nicer to compute, to maximize
$$\left(1-\frac{k}{a_1}\right) a_{k+1}$$
   under the condition
$$ka_1^{N-1}+a_{k+1}^{N-1} \leq r.$$
   Since we want to go down with only one case, we force $a_{k+1}$ to be an integer. Now the problem is to distribute
points to $r_1$ and $r_{k+1}$. It is a matter of an easy calculation to check integer $a_{k+1}$ with $r_1=\lfloor (r-a_{k+1}^{N-1})/k \rfloor$ gives the best result.

In our case, $N=4$ and $r=180$, the following distribution was found:
$$r_1 = 51, \qquad r_2 = 51, \qquad r_3 = 51, \qquad r_4 = 27.$$
   Thus we need a lower bound for $\alphahat(\P^3;51)$. Again, we use the above heuristic method to find the distribution
   $$r_1' = 14, \qquad r_2' = 14, \qquad r_3' = 14, \qquad r_4' = 9.$$
   We take the bound for $\alphahat(\P^2;14) \geq 86/23$. Thus
$$\alphahat(\P^3;51) \geq \frac{309}{86}, \qquad \alphahat(\P^4;180) \geq \frac{360}{103} \simeq 3.495.$$
   Our previous best bound is better only by $\simeq 0.0003549$ but the run time of the algorithm outlined here is considerably shorter.

Less radical, but a better approach is to consider all distributions $kr_1+r_{k+1} \leq r$ with $r_{k+1}$ being a pure $(N-1)$th power. The implementation of these two approaches in Singular \cite{Singular} can be found in the file \verb"boundforWC", \cite{boundforWC}.
Running \verb"bound" works faster (for big $N$), but \verb"boundmore" gives better bounds.
\end{example}

\subsection{An easy way to distribute points on hyperplanes}
   We pass now to some general effective lower bounds.
\begin{proposition}\label{prop:distribution k s-k}
   Let $k$ be a positive integer and let $s$ be an integer in the range $1\leq s\leq k$. Let
   $$r \geq s(k+1)^{N-1}+(k+1-s)k^{N-1}.$$
   Then
   $$\alphahat(\P^N;r) \geq k+\frac{s}{k+1}.$$
\end{proposition}
\proof
   This is an easy consequence of Theorem \ref{stepback}. Namely,
taking
$$r_1 = \ldots = r_s = (k+1)^{N-1}, \qquad r_{s+1}= \ldots = r_{k+1} = k^{N-1},$$
we get by Theorem \ref{thm:lower bound}
$$a_1 = \ldots = a_s = k+1, \qquad a_{s+1} = \ldots = a_{k+1} = k.$$
Consequently,
$$\alphahat(r) \geq \left(1-\frac{s}{k+1}-\frac{k-s}{k}\right)k+k = s-\frac{sk}{k+1}+k=k+\frac{s}{k+1}.$$
\endproof

\begin{example}
Without the above proposition, the general available lower bound for $\alphahat(\P^5;1024)$ is $4$.
   It requires at least $r \geq 3125$ points to pass to the better bound $\alphahat(\P^5;r) \geq 5$. But
   with Proposition \ref{prop:distribution k s-k} we can take
$s=1$, $k=4$ to get
$$\alphahat(\P^5;1649) \geq 4 + \frac{1}{5}.$$
Similarly, we need only $2018$ points to get $4 + \frac{2}{5}$, only $2387$ to get $4 + \frac{3}{5}$ and only $2756$
to get $4 + \frac{4}{5}$.
\end{example}

\begin{proposition}\label{quotb}
Let $r \leq (k+1)^N$. Then
$$\alphahat(\P^N;r) \geq \frac{r}{(k+1)^{N-1}}.$$
\end{proposition}

\proof
We will use induction. Consider two cases.

\medskip

\textbf{Case $r \leq k(k+1)^{N-1}$.}\\
Since
$$k \leq k(k+1)^{N-1} \leq k^N,$$
we have (by induction on $k$)
$$\alphahat(\P^N;r) \geq \frac{r}{k^{N-1}} \geq \frac{r}{(k+1)^{N-1}}.$$

\medskip

\textbf{Case $r > k(k+1)^{N-1}$.}\\
Take
$$r_1 = \ldots = r_k = (k+1)^{N-1}, \qquad r_{k+1} = r-k(k+1)^{N-1}.$$
Observe that
$$r-k(k+1)^{N-1} \leq (k+1)^N-k(k+1)^{N-1}=(k+1)^{N-1},$$
thus by Theorem \ref{thm:lower bound} and induction (on $N$) we get
$$a_1 = \ldots = a_k = k+1, \qquad a_{k+1} = \frac{r-k(k+1)^{N-1}}{(k+1)^{N-2}}.$$
By Theorem \ref{stepback} we get
$$\alphahat(\P^N;r) \geq \left(1-\frac{k}{k+1}\right)\frac{r-k(k+1)^{N-1}}{(k+1)^{N-2}}+k = \frac{r}{(k+1)^{N-1}}.$$
\endproof

\begin{example}
We can now give very accurate bounds for $\alphahat(\P^5;r)$ for $r$ close to $3125$.
Since $3125 = 5^5$, we have
$$\alphahat(\P^5;3124) \geq 5 - \frac{1}{625}, \qquad \alphahat(\P^5;3123) \geq 5 - \frac{2}{625}$$
and so on.
\end{example}

\subsection{Discussion on the accuracy}

\begin{figure}[h]
\caption{Upper and lower bounds for $\alphahat(\P^3;r)$}
\label{figp3}
\centering
\includegraphics[scale=0.3]{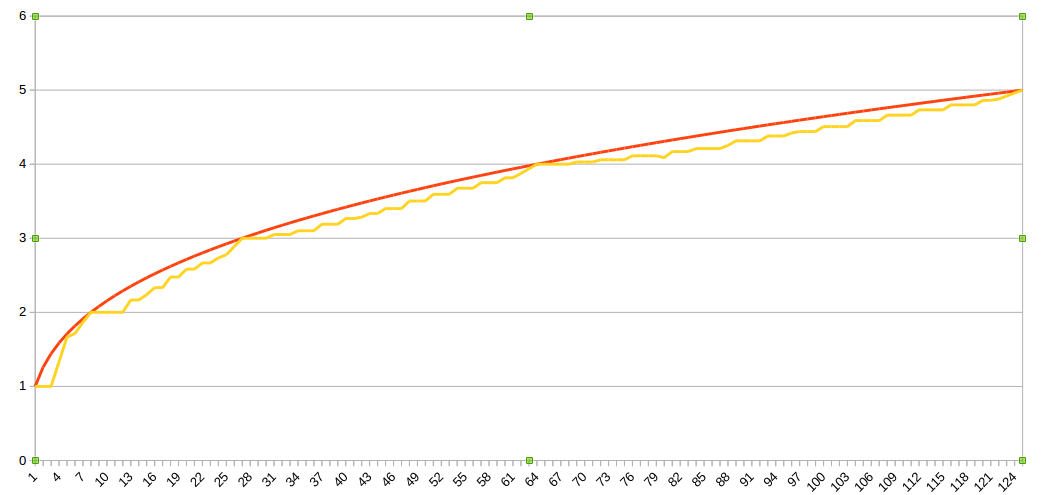}
\end{figure}

By Theorem \ref{thm:lower bound} it is obvious that we can locate every $\alphahat(\P^N;r)$ in an interval of
length at most $1$. It is interesting to know what is the difference between the upper bound (which is conjectured
to be the actual bound for $r \geq 2^N$) and the lower bound obtained by our algorithm. In Figure \ref{figp3}
we present the upper and lower bounds for $r=1,\dots,125$ points in $\P^3$.

In Table \ref{tab:accuracy} we present the maximal difference $\delta$ between the lower and upper bound.
\renewcommand*{\arraystretch}{1.2}
\begin{table}[h]
$$
\begin{array}{c|cc|cc|cc|c|c}
N  &      3    &     3     &    4     &    4      &     5      &      5     &     6     &     7     \\ \hline
r  &  8-125   & 125-1000 &  16-256 & 256-1296 & 32-243    & 243-1024  & 64-729   & 128-2187 \\
\delta & 0.289 &  0.186    &  0.295   & 0.259     &   0.305    & 0.277      &  0.305    &  0.301    \\
\end{array}
$$
\caption{Maximal differences for lower and upper bounds in various intervals of the number of points in projective spaces of low dimensions}
\label{tab:accuracy}
\end{table}

\subsection{Towards Demailly's Conjecture}
   As an important consequence of Theorem \ref{thm:main} we obtain the following result.
\begin{theorem}\label{thm:DC ok for geq mN}
   Demailly's Conjecture \ref{conj:Demailly} holds for $r\geq m^N$ very general points in $\P^N$.
\end{theorem}
\proof
   The Main Theorem in \cite{MSS17} states that Conjecture \ref{conj:Demailly}
   holds for $r\geq (m+1)^N$ very general points in $\P^N$. Hence it is enough
   to deal with sets $Z$ containing $r$ very general points with $r$ in the range $m^N\leq r<(m+1)^N$.
   The general yoga of our proof is the following: We use lower bounds on the Waldschmidt constant of $Z$
   provided either by Theorem \ref{thm:lower bound} or by Proposition \ref{prop:distribution k s-k}
   and check, by naive conditions count, that $\alpha(mZ)$ is small enough for the inequality
   \eqref{eq:Chudnovsky Conjecture} to be satisfied.\\
   \textbf{Case 1.} For $N\geq 4$ and $m\geq 3$, it follows from Lemma \ref{lem:inequality}
   that there exists a hypersurface in $\P^N$ of degree $m(m+N-1)-N+1$
   vanishing to order at least $m$ at all points of $Z$. Since in any case it is $\alphahat(Z)\geq m$
   by Theorem \ref{thm:lower bound}, it follows that
   $$\alphahat(Z)\geq m\geq \frac{\alpha(mZ)+N-1}{m+N-1}$$
   and we are done in this case.\\
   \textbf{Case 2.} Let $N=3$ and let $2\leq m=2n+\eps$ with $\eps\in\left\{0,1\right\}$.
   Assume that
   $$m^3\leq r\leq (n+1+\eps)(m+1)^2+(n-1)m^2.$$
   It follows again from the naive conditions count that there exists a surface in $\P^3$
   of degree $m^2+2m-2$ passing with multiplicity at least $m$ through all points in $Z$.
   Hence $\alpha(mZ)\leq m(m+2)-2$ and thus
   $$\frac{\alpha(mZ)+2}{m+2}\leq m\leq \alphahat(Z),$$
   which is exactly \eqref{eq:Demailly Conjecture}.\\
   If the number of points is in the range
   $$(n+1+\eps)(m+1)^2+(n-1)m^2\leq r <(m+1)^3,$$
   then Proposition \ref{prop:distribution k s-k} implies that
   $$\alphahat(Z)\geq m+1-\frac{m+1-n-1-\eps}{m+1}\geq m+\frac12.$$
   If $m=2n$ is even, then there exists a surface of degree $4n^2+9n+2$
   vanishing at all points of $Z$ to order at least $m$. Indeed, this follows
   from the inequality
   $$\binom{4n^2+9n+5}{3}\geq (m+1)^3\binom{m+2}{3},$$
   which is equivalent to
   $$8n^5+(62/3)n^4+(39/2)n^3+(47/6)n^2+n\geq 0.$$
   Hence $\alpha(mZ)\leq (m+\frac12)(m+2)-2$, which gives
   $$\alphahat(Z)\geq m+\frac{1}{2}\geq \frac{\alpha(mZ)+2}{m+2},$$
   hence \eqref{eq:Demailly Conjecture} holds.\\
   The case $m=2n+\eps$ is similar and we leave it as a simple exercise.\\
   \textbf{Case 3.} Let $m=2$ and let $Z$ be a set of $r$ very general points
   in $\P^N$ with $2^N\leq r<3^N$. In any case it is $\alphahat(Z)\geq 2$
   by Theorem \ref{thm:lower bound}. For $N\geq 7$ this bound suffices
   to conclude Conjecture \ref{conj:Demailly}. Indeed, since
   $$\binom{2N+3}{N}\geq 3^N(N+1)\;\;\mbox{ holds for }N\geq 7,$$
   there is a hypersurface of degree $N+3$ singular in points of $Z$.
   Hence $\alpha(2Z)\leq N+3$ and this implies
   $$\alphahat(Z)\geq 2\geq\frac{\alpha(2Z)+N-1}{N+1}.$$
   For $4\leq m\leq 6$ we split the argument in two cases:
   \begin{itemize}
   \item[a)] $r\leq 2\cdot 3^{N-1}+2^{N-1}$ and
   \item[b)] $r\geq 2\cdot 3^{N-1}+2^{N-1}$.
   \end{itemize}
   In case a) the previous argument works. There is a hypersurface
   of degree $N+3$ in $\P^N$ singular in points of $Z$.
   In case b) we
   apply Proposition \ref{prop:distribution k s-k}
   with $s=2$ and $k=2$. It follows then that $\alphahat(Z)\geq 8/3$.
   By elementary conditions count, there is a hypersurface of degree $2N+1$
   singular at $Z$, so that $\alpha(2Z)\leq 2N+1$. Hence
   $$\alphahat(Z)\geq \frac83\geq\frac{2N+1+N-1}{N+1}$$
   holds as $N\leq 7$. \\
   \textbf{Case 4.} Finally we are left with $m=1$ but this has been proved
   for all $N$ in \cite{DT16} and independently in \cite{FMX16}.
\endproof

\begin{remark}
Using similar methods one can easily check if the bound for $\alphahat(\P^N;r)$ is
sufficient to prove the Demailly Conjecture for a given $N$, $m$ and $r$. We wrote
an appropriate procedure (\verb"Demailly" in \verb"boundforWC") and check that, for example,
the Conjecture holds for all $N \leq 3$, $m \leq 3$ and any number of very general points.
\end{remark}

\section{Auxiliary results}

\begin{lemma}\label{solutionlem}
Assume that positive real numbers $a_1,\dots,a_{t-1}$ are given, satisfying
$$1-\sum\limits_{j=1}^{t-1} \frac{1}{a_j} \neq 0.$$
Let $C$ be a real number. Consider the following system of linear equations:
$$\left\{
\begin{array}{rcl}
C - \sum\limits_{i=1}^{t-1} x_i & = & a_k(1-x_k) \text{ for } k=1,\dots,t-1 \\
y & = & \sum\limits_{i=1}^{t-1} x_i.
\end{array}
\right.$$
   Then there is the unique solution for $(x_1,\dots,x_{t-1},y)$ to this system. In particular
$$y = \frac{C(\sum\limits_{j=1}^{t-1} \frac{1}{a_j})-(t-1)}{\sum\limits_{j=1}^{t-1} \frac{1}{a_j} - 1}.$$
\end{lemma}

\proof
We look for the (unique) solution for $y$, thus we use Cramer's rule. The matrix of this system (after some reorganisation:
the variable $y$ is placed in the first column, then $x_1,\dots,x_{t-1}$, then non-linear part) is equal to
$$M := \left[\begin{array}{ccccccc}
0  &  1-a_1  &  1  &  1 &   \dots &  1 & C-a_1 \\
0  &   1    & 1-a_2 & 1 &  \dots &  1 & C-a_2 \\
0  &   1    &  1  & 1-a_3 & \dots & 1 & C-a_3 \\
\vdots & \vdots & \vdots & \vdots & \vdots & \vdots & \vdots \\
0  &   1    &  1  & 1   & \dots & 1-a_{t-1} & C-a_{t-1} \\
-1 &   1    &  1  & 1   & \dots &  1    &  0.
\end{array}
\right].$$

We denote the columns of $M$ by $[A_0 A_1 \dots A_{t-1} B]$.

To compute the determinant of the main matrix $[A_0 A_1 \dots A_{t-1}]$ we subtract the last row from the others, obtaining the matrix
with the first column and last rows filled with $1$ (except $-1$ in the left bottom corner), and then $-a_1,-a_2,\dots$ over the diagonal.
Applying Laplace rule we compute this determinant to be equal to
$$D_1 = \left( \sum_{i=1}^{t-1} a_1\ldots \widehat{a_i} \ldots a_{t-1} \right) - a_1 \ldots a_{t-1} = a_1 \ldots a_{t-1} \cdot \left(\sum\limits_{j=1}^{t-1} \frac{1}{a_j} - 1 \right)$$
which is non-zero (by the assumption). Hence the solution is unique.

To compute the determinant of the matrix $[B A_1 A_2 \dots A_{t-1}]$ we ''kill'' all $1$'s using the last row, then ''kill''
all $a_i$'s in the first column using other columns, obtaining the matrix with
$$\left[\begin{array}{cccc}
C & -a_1 & & \\
\vdots & & \ddots & \\
C & & & -a_{t-1} \\
t & 1 & \ldots & 1
\end{array}
\right].$$
By the Laplace rule, the determinant
$$D_2 = C \left( \sum\limits_{i=1}^{t-1} a_1 \ldots \widehat{a_i} \ldots a_{t-1} \right) - (t-1) a_1 \dots a_{t-1} =
a_1 \dots a_{t-1} \left( C \left( \sum_{j=1}^{t-1} \frac{1}{a_j} \right) - (t-1) \right).$$
By the Cramer's rule, the claim follows.
\endproof

\begin{lemma}\label{lem:inequality}
   For all $N\geq 4$, $m\geq 3$ there is
   \begin{equation}\label{eq: combinatorial}
      \binom{m(m+N-1)+1}{N}>\binom{m+N-1}{N}(m+1)^N.
   \end{equation}
\end{lemma}
\proof
   With $m\geq 3$ fixed, the proof goes by induction on $N$. In the initial case $N=4$ it is elementary
   to check that the claim is equivalent to the inequality
   $$m^2(2m^5+11m^4-89m^2-146m-42)>0,$$
   which is fulfilled for all $m\geq 3$. \\
   For the induction step, we assume that \eqref{eq: combinatorial} holds and we want to
   show that
   \begin{equation}\label{eq:loc4}
      \binom{m(m+N)+1}{N+1}> \binom{m+N}{N+1}(m+1)^{N+1}
   \end{equation}
   holds as well. It is convenient to abbreviate $A=m(m+N)$. Using the induction assumption
   and after elementary operations we get
   \begin{align*}
      \binom{m(m+N)+1}{N+1}&> \binom{m+N}{N+1}(m+1)^{N+1}\cdot\\
         &\cdot\frac{(A+1)A(A-1)\ldots(A-m)}{(m+1)(m+N)(A-N)(A-N-1)\ldots(A-N-m+2)},
   \end{align*}
   so that in order to get \eqref{eq:loc4}, it suffices to show
   $$(A+1)A(A-1)\ldots(A-m)\geq (m+1)(m+N)(A-N)(A-N-1)\ldots(A-N-m+2),$$
   which follows by comparing both sides term by term (there are $(m+1)$ terms
   on both sides of the inequality).
\endproof
\paragraph*{Acknowledgement.}
   Our research was partially supported by National Science Centre, Poland, grant
   2014/15/B/ST1/02197.


\bigskip \small

\bigskip
   Marcin, Dumnicki,
   Jagiellonian University, Institute of Mathematics, {\L}ojasiewicza 6, PL-30348 Krak\'ow, Poland

\nopagebreak
   \textit{E-mail address:} \texttt{marcin.dumnicki@im.uj.edu.pl}

\bigskip

   Tomasz Szemberg,
   Department of Mathematics, Pedagogical University of Cracow,
   Podchor\c a\.zych 2,
   PL-30-084 Krak\'ow, Poland.

\nopagebreak
   \textit{E-mail address:} \texttt{tomasz.szemberg@gmail.com}

\bigskip

   Justyna Szpond,
   Department of Mathematics, Pedagogical University of Cracow,
   Podchor\c a\.zych 2,
   PL-30-084 Krak\'ow, Poland.

\nopagebreak
   \textit{E-mail address:} \texttt{szpond@gmail.com}


\end{document}